\documentclass[12pt,reqno]{amsart}

\usepackage{amssymb,amsmath, amsfonts, latexsym}
\usepackage{enumerate}
\newtheorem{thm}{Theorem}[section]
\newtheorem{cor}[thm]{Corollary}
\newtheorem{lem}[thm]{Lemma}
\newtheorem{prop}[thm]{Proposition}
\newtheorem{example}[thm]{Example}
\newtheorem{remarks}[thm]{Remark}
\newtheorem{defn}[thm]{Definition}

 \numberwithin{equation}{section}

\title[Bernstein type's concentration inequalities]
{Bernstein type's concentration inequalities for symmetric Markov processes}

\author{Fuqing Gao}
\address{Fuqing GAO.
School of Mathematics and Statistics, Wuhan University, 430072
Hubei, China} \email{fqgao@whu.edu.cn}

\author{Arnaud Guillin}
\address{Arnaud Guillin. Laboratoire de Math\'ematiques Appliqu\'ees, CNRS-UMR
6620, Universit\'e Blaise Pascal, 63177 Aubi\`ere, France. }
\email{Arnaud.Guillin@math.univ-bpclermont.fr}

\author{Liming Wu}
\address{Liming Wu. Laboratoire de Math\'ematiques Appliqu\'ees, CNRS-UMR
6620, Universit\'e Blaise Pascal, 63177 Aubi\`ere, France. And
Institute of Applied Mathematics, Chinese Academy of Sciences,
100190 Beijing, China.} \email{Li-Ming.Wu@math.univ-bpclermont.fr}

\date{May 2009}

\newcommand{\dd}{\mathbb{D}}
\newcommand{\ee}{\mathbb{E}}

\newcommand{\nn}{\mathbb{N}}
\newcommand{\rr}{\mathbb{R}}
\newcommand{\pp}{\mathbb{P}}
\newcommand{\qq}{\mathbb{Q}}

\newcommand{\Var}{\mathrm{Var}}

\def\BB{\mathcal B}
\def\CC{\mathcal C}

\def\FF{\mathcal F}

\def\EE{\mathcal E}

\def\LL{\mathcal L}
\def\MM{\mathcal M}
\def\NN{\mathcal N}

\def\XX{\mathcal X}

\def\vep{\varepsilon}
\def\<{\langle}
\def\>{\rangle}

\def\beq{\begin{equation}}
\def\neq{\end{equation}}

\def\bthm{\begin{thm}}
\def\nthm{\end{thm}}

\def\bprop{\begin{prop}}
\def\nprop{\end{prop}}

\def\brmk{\begin{remarks}}
\def\nrmk{\end{remarks}}

\def\bexa{\begin{example}}
\def\nexa{\end{example}}

\def\blem{\begin{lem}}
\def\nlem{\end{lem}}

\def\bcor{\begin{cor}}
\def\ncor{\end{cor}}

\def\bexe{\begin{exe}}
\def\nexe{\end{exe}}

\def\bprf{\begin{proof}}
\def\nprf{\end{proof}}

\def\bdef{\begin{defn}}
\def\ndef{\end{defn}}

\def\dsp{\displaystyle}

\def\bdes{\begin{description}}
\def\ndes{\end{description}}

\def\benu{\begin{enumerate}}
\def\nenu{\end{enumerate}}

\def\Var{\text{\rm Var}}
\def\Ent{\text{\rm Ent}}

\begin{document}

\begin{abstract} Using the method of transportation-information
inequality introduced in \cite{GLWY}, we establish Bernstein type's
concentration inequalities for empirical means $\frac 1t \int_0^t
g(X_s)ds$ where $g$ is a unbounded observable of the symmetric
Markov process $(X_t)$. Three approaches are proposed : functional
inequalities approach ; Lyapunov function method ; and an approach
through the Lipschitzian norm of the solution to the Poisson
equation. Several applications and examples are studied.
\end{abstract}

\maketitle




\section{Introduction}

\subsection{Bernstein's concentration inequality for sequences of
i.i.d.r.v. } Let us begin with the classical Bernstein's
concentration inequality in the i.i.d. case. Consider a sequence of
real valued independent and identically distributed (i.i.d.) random
variables (r.v.) $(\xi_k)_{k\ge 1}$, copies of some r.v. $\xi$, all
defined on the probability space $(\Omega,\FF,\pp)$ such that $\ee
\xi=0$ and $\ee \xi^2=\sigma^2>0$.

\bthm\label{thm-Bernstein} If there is some constant $M\ge0$ such
that \beq\label{Bern-a} \Lambda(\lambda):=\log \ee e^{\lambda \xi}
\le \frac{\lambda^2 \sigma^2}{2(1-\lambda M)}, \ \lambda\in (0,1/M).
\neq Then for any $r>0$ and $n\ge 1$,
 \beq\label{Bern-b} \pp\left(\frac 1n
\sum_{k=1}^n \xi_k
> r\right)\le \exp\left(-n
\frac{2r^2}{\sigma^2\left(\sqrt{1+\frac{2Mr}{\sigma^2}}+1\right)^2}\right),
\ r>0 \neq or equivalently for any $x>0$ and $n\ge 1$,
\beq\label{Bern-c} \pp\left(\frac 1n \sum_{k=1}^n \xi_k > \sigma
\sqrt{2x}+ Mx\right)\le e^{-n x}. \neq  In particular
\beq\label{Bern-d} \pp\left(\frac 1n \sum_{k=1}^n \xi_k >
r\right)\le \exp\left(- \frac{nr^2}{2(\sigma^2 +Mr)}\right), \ r>0.
\neq \nthm

The last inequality (\ref{Bern-d}) is the original version of
Bernstein's inequality. The proof of (\ref{Bern-b}) is very easy :
just apply Chebychev's inequality to obtain : $\forall r,\lambda>0$,
$$
\pp\left(\frac 1n \sum_{k=1}^n \xi_k
> r\right)\le e^{-n\lambda r} \ee \exp\left(\lambda \sum_{k=1}^n \xi_k\right)\le e^{-n[\lambda r - \Lambda(\lambda)]}
$$
and then optimize over $\lambda\in (0, 1/M)$. We refer to E. Rio
\cite{Rio00} or P. Massart \cite{massart} for known sufficient
conditions for the verification of (\ref{Bern-a}). For instance
(\ref{Bern-a}) is verified with $M=\|\xi^+\|_\infty/3$ if $\xi$ is
upper bounded, or for some not very explicit  constant  $M>0$ if
$\Lambda(\lambda)<+\infty$ for {\bf some} $\lambda>0$. Bernstein's
concentration inequality is one of the most powerful concentration
inequalities in probability, which is sharp both in the central
limit theorem scale and the moderate deviation scale. This type of
inequalities have had many applications, and are now particularly
used in (non asymptotic) model selection problem, see Massart
\cite{massart} or Baraud \cite{Baraud}.

There are already many works on the generalization of Bernstein's
inequality in the dependent case: Markov process or weakly dependent
one. The strategy however remains the same : control the Laplace
transform of partial sums. In the markovian context, Lezaud
\cite{Lez01} used Kato's perturbation theory to get result in
presence of a spectral gap, whereas Cattiaux-Guillin \cite{CG08}
(building on Wu \cite{Wu00c}) used functional inequalities for the
Laplace control or for the control of the mixing coefficients. More
recently, Adamczak \cite{Ad08}, Bertail-Cl\'emen\c con \cite{BC08},
Merlev\`ede-Peligrad-Rio \cite{MPR09} used a block strategy and then
results in the independent case. Note however that, except the
symmetric Markov processes case studied by Lezeaud \cite{Lez01}, the
known results do not reach the tight form (\ref{Bern-b}) or
(\ref{Bern-d}). Our major objective is to give practical conditions
ensuring this sharp form (\ref{Bern-b}) in the context of integral
functional of symmetric Markov processes.

There are two modern approaches to concentration inequalities. The
first one, initiated by Ledoux,  relies on functional inequalities,
such as Poincar\'e or logarithmic Sobolev inequality (see for
example \cite{logsob} or \cite{Led01}) and has attracted a lot of
attention in the past decade: Wu \cite{Wu00c} or Cattiaux-Guillin
\cite{CG08} used them in the continuous time context to get precise
control of the Laplace transform of the partial sums, see also
Massart \cite{massart} for the entropy method for various type of
dependance in the discrete time case; another approach was to get a
functional inequality for the whole law of the process and Herbst's
like argument, note however that at this level of generality, the
precise form of Bernstein's inequality has not been achieved yet.

The second approach is centered on the use of transportation
inequalities ( see precise definition in section 2 below): bounding
Wasserstein's distance by some type of information (Kullback or
Fisher). If originally investigated by Marton \cite{Mar96,Mar97} or
Talagrand \cite{Tal96a} for concentration, its systematic study is
more recent, starting from the pioneer work of Bobkov-Gotze
\cite{BG99}, followed by an abundant litterature, see
\cite{OVill00,BGL01,DGW03,BV05,CG06,GL} with Kullback information,
and \cite{GLWY,GLWW,GJLW} for Fisher information. If the use of
Kullback information at the process level may lead to deviation
inequality for integral functional of Markov processes (see
\cite{DGW03} for example), the precise form of Bernstein's
inequality is not reachable. We will therefore use here
transportation inequalities with respect to the Fisher information,
which are more natural for Markov processes : the Fisher information
is exactly the large deviations rate in the Donsker-Varadhan theorem
for symmetric Markov processes (see
\cite{DV75a,DV76,DV83,Wu00b,KM}).

But before going further into the details, let us present the framework on symmetric Markov processes.


\subsection{Symmetric Markov processes}

Let $\XX$ be a Polish space with Borel field $\BB$. Let $(X_t)_{t\ge
0}$ be a $\XX$-valued {\it c\`adl\`ag} Markov process with
transition probability semigroup $(P_t)$ which is symmetric and
strongly continuous on $L^2(\mu):=L^2(E,\BB,\mu)$, defined on
$(\Omega, \FF, (\pp_x)_{x\in \XX})$ $(\pp_x(X_0=x)=1,\ \forall
x\in\XX)$, where $\mu$ is a probability measure on $(\XX,\BB)$,
written as $\mu\in \MM_1(\XX)$. For a given initial distribution
$\beta\in\MM_1(\XX)$, write $\pp_\beta:=\int_\XX \beta(dx)
\pp_x(\cdot)$. Let $\LL$ be the generator of $(P_t)$, whose domain
in $L^p(\mu)=L^p(\XX,\BB,\mu)\ (p\in [1,+\infty])$ is denoted by
$\dd_p(\LL)$. It is self-adjoint, definitely non-positive on
$L^2(\mu)$. Let
$$
-\LL=\int_0^{+\infty} \lambda dE_\lambda
$$
be the spectral decomposition of $-\LL$ on $L^2(\mu)$. The Dirichlet
form $\EE(f,g)$ is defined by
$$\aligned
&\dd(\EE)=\dd_2(\sqrt{-\LL})=\left\{h\in L^2(\mu);\ \int_0^{+\infty}
\lambda d\<E_\lambda h, h\>_\mu <+\infty\right\}\\
&\EE(f,g)=\<\sqrt{-\LL}f, \sqrt{-\LL}g\>_\mu=\int_0^{+\infty}
\lambda d\<E_\lambda f, g\>_\mu,\ f,g\in\dd(\EE)
 \endaligned$$
where $\<f,g\>_\mu=\int_\XX fg d\mu$ is the standard inner product
on $L^2(\mu)$.

\smallskip

We will study here deviation inequalities for
$$\frac1t\int_0^tg(X_s)ds$$ for some $\mu$-centered function $g$ (observable).
It is quite natural to expect conditions relying on an interplay
between the type of ergodicity of our Markov process and the type of
boundedness or integrability of the function $g$.

That is why a long standing assumption in this paper will be the following Poincar\'e inequality
: for some finite nonnegative best constant $c_P$,
\beq\label{Poincare} \Var_\mu(f)\le c_P \EE(f,f),\ \forall f\in
\dd(\EE). \neq Here and hereafter $\mu(f):=\int_\XX f d\mu$ and
$\Var_\mu(f)=\mu(f^2)-\mu(f)^2$ is the variance of $f$ under $\mu$.
Poincar\'e's inequality is equivalent to the exponential decay of
$P_t$ to the equilibrium invariant measure $\mu$ in $L^2(\mu)$ :
$$
\Var_\mu(P_t f)\le e^{-2t/c_P} \Var_\mu(f),\ \forall f\in L^2(\mu).
$$
It is also equivalent to say that the spectral gap
$$
\lambda_1:=\sup\{\lambda\ge 0; \ E_\lambda-E_0=0\}=\frac 1{c_P}>0.
$$

Let us first show why this Poincar\'e inequality condition is
natural in our context. Indeed, the first class of test function $g$
that can be considered is the class of bounded ones. Using Kato's
theory about perturbation of operators combined with ingenious and
difficult combinatory calculus, Lezaud \cite{Lez01} proved the
following Bernstein type's concentration inequality.

\bthm\label{thm-Lezaud}{\rm (\cite{Lez01})} Let $g$ be a  bounded
and measurable function (say $g\in b\BB$) such that $\mu(g)=0$. Then
for $\beta\ll \mu$,

\beq\label{Bern1}\aligned \pp_\beta\left(\frac 1t \int_0^t
g(X_s)ds>r\right)&\le \|\frac{d\beta}{d\mu}\|_2 \exp\left(-
\frac{2tr^2}{\sigma^2\left(\sqrt{1+\frac{2Mr}{\sigma^2}}+1\right)^2}\right)\\
&\le \|\frac{d\beta}{d\mu}\|_2\exp\left(- \frac{tr^2}{2(\sigma^2
+Mr)}\right), \ \forall \ t, r>0\endaligned\neq where $M=M(g)=c_P
\|g\|_\infty$ and $\sigma^2$ is  the asymptotic variance (in the
CLT) of the observable $g\in L^2(\mu)$,  given by
 \beq\label{variance}
\sigma^2=\sigma^2(g):=\lim_{t\to+\infty} \frac 1t \Var_{\pp_\mu}
\left(\int_0^t g(X_s)ds\right)= 2 \int_0^{+\infty} \<P_t g, g\>_\mu
dt. \neq
 \nthm

 For generalization of this result see
Cattiaux-Guillin \cite{CG08}, Guillin-L\'eonard-Wu-Yao \cite{GLWY}
etc. Notice a remarkable point : (\ref{Bern1}) is sharp both for the
central limit theorem (CLT) scale $r\propto 1/\sqrt{t}$ (since
$\frac 1{\sqrt{t}}\int_0^t g(X_s)ds$ converges in law to the
centered Gaussian distribution with variance $\sigma^2(g)$, see \cite{KV86}), and for
the moderate deviation scale (i.e. $1/\sqrt{t}\ll r\ll 1$) by the
moderate deviation principle due to \cite{Wu95}.

Notice that if $\sigma^2(g)\le C \|g\|_\infty^2$ for some constant
$C>0$ and for all $g\in b\BB$ with $\mu(g)=0$, then the Bernstein's
concentration inequality (\ref{Bern1}) implies the Poincar\'e
inequality (\ref{Poincare}), by \cite[Theorem 3.1]{GLWY}. In other
words the Poincar\'e inequality is a minimal assumption for
Bernstein's concentration inequality for {\bf all} bounded
observables $g$.

\brmk{\rm
Let us point out that for bounded $g$, the assumption that $\sigma^2(g)\le C \|g\|_\infty^2$ is a weak one, as by definition (\ref{variance})
$$\sigma^2(g)\le 4\|g\|_\infty\int_0^t\Var_\mu(P_tg)^{1/2}dt.$$
Assume now that a weak Poincar\'e inequality holds (see \cite{BCG2} for example), or a Lyapunov condition, i.e. $\LL V\le-\phi(V)+b1_C$ for some sub linear $\phi$ (see \cite{DFG} for details), ensuring that $\Var_\mu(P_tg)\le \psi(t)\|g\|^2_\infty$ with $\int_0^s\psi(s)^{1/2}ds<\infty$, then the Poincar\'e inequality holds under Bernstein's type inequality. We refer to the last section for some examples of this Lyapunov condition.
}\nrmk


\subsection{Main question and organization}
The main question we will focus on in this paper will be:
 {\it what is the interplay between the ergodic properties of the symmetric Markov process and the test function $g$?
}
 Or more precisely, how to bound the constant $M$ (appearing in (\ref{Bern1})) by means of other quantities
 than $\|g\|_\infty$ and $c_P$?

 In fact we shall answer this question by a very simple approach : instead of a direct
control of the Laplace transform of partial sums, we use the method
of  transportation-information inequality introduced by
Guillin-L\'eonard-Wu-Yao \cite{GLWY}.

This paper is organized as follows. In the next section we describe
the strategy and the main idea of this work, giving by the way
another proof of Theorem \ref{thm-Lezaud} with a better estimate of
$M$. The goal of the three following sections is to generalize
Bernstein's inequality to unbounded case. We present three
approaches : (1) functional inequalities such as log-Sobolev
inequality or $\Phi$-Sobolev inequality ; (2) the Lipschitzian norm
$\|(-\LL)^{-1}g\|_{Lip}$ ; and (3) Meyn-Tweedie's Lyapunov function
method. Finally the last section is dedicated to the case where
Poincar\'e inequality does not hold anymore, and the class of
bounded test functions is now too large. Once again, the approach
via Lyapunov function will be particularly efficient.

Note that, from {\it Section 2 through 5}, we assume implicitly that the previous Poincar\'e inequality is satisfied.

Before going to the job let us fix some more notations. For $p\in
[1,+\infty]$, $\|\cdot\|_p$ is the standard norm of
$L^p(\mu):=L^p(\XX,\BB,\mu)$, and $L^p_0(\mu):=\{g\in L^p(\mu);\
\mu(g)=0\}$. The quantity $\sigma^2$ denotes always the asymptotic
variance $\sigma^2(g)$ in the CLT,  given by (\ref{variance}). The
empirical measure $\frac 1t \int_0^t \delta_{X_s} ds$ ($\delta_x$
being the Dirac measure at point $x$) is denoted by $L_t$, so that
$\frac 1t \int_0^t g(X_s) ds =L_t(g)$.


\section{A transportation-information look at Bernstein's inequality}
\subsection{The strategy and the main idea} As in \cite{GLWY}, our starting point is

\bthm\label{prop21} {\rm (Wu \cite{Wu00c})} Let $g\in L^1_0(\mu)$.
Then \beq\label{prop21a} \pp_\beta\left(\frac 1t \int_0^t
g(X_s)ds>r\right)\le \|\frac{d\beta}{d\mu}\|_2e^{-t I(r-)},\ \forall
t,r>0 \neq where
$$
I(r):=\inf\{I(\nu|\mu);\ \nu(|g|)<+\infty, \nu(g)=r\}, \
I(r-):=\lim_{\vep\to 0+} I(r-\vep), \ r\in \rr
$$
and \beq\label{Fisher} I(\nu|\mu):=\begin{cases} \EE\left(\sqrt{f},
\sqrt{f}\right), \ &\text{ if }\ \nu=f\mu, \sqrt{f}\in \dd(\EE),\\
+\infty, &\text{ otherwise}
\end{cases} \neq
is the Fisher-Donsker-Varadhan's information of $\nu$ with respect
to (w.r.t.) $\mu$. \nthm

By the large deviations in Donsker-Varadhan \cite{DV75a, DV76} (in
the regular case) and Wu \cite{Wu00b} (in full generality), $\nu\to
I(\nu|\mu)$ is the rate function in the large deviations of the
empirical measures $L_t:=\frac 1t \int_0^t \delta_{X_s} ds$, and the
Cramer type's inequality (\ref{prop21a}) is sharp for large time
$t$. The main problem now is to estimate the rate function $I(r)$ in
the large deviations of $\frac 1t \int_0^t g(X_s)ds$ : that is
exactly a role that the transportation-information inequality plays.

\bthm\label{thm22} {\rm (\cite[Theorem 2.4]{GLWY})} Let $g\in
L^1_0(\mu)$ and $\alpha: \rr\to [0,+\infty]$ be a nondecreasing
left-continuous convex function with $\alpha(0)=0$. The following
properties are equivalent :

\begin{enumerate}[(a)]
\item $\alpha(\nu(g))\le I(\nu|\mu),\ \forall \nu\in\MM_1(\XX)$ such that $\nu(|g|)<+\infty$.
\item $\nu(g)\le \alpha^{-1}(I(\nu|\mu)),\ \forall \nu\in\MM_1(\XX)$ such that $\nu(|g|)<+\infty$,
where $\alpha^{-1}(x):=\inf\{r\in \rr;\ \alpha(r)>x\}$ is the right
inverse of $\alpha$.
\item It holds that
\beq\label{thm22a} \pp_\beta\left(\frac 1t \int_0^t
g(X_s)ds>r\right)\le \|\frac{d\beta}{d\mu}\|_2e^{-t \alpha(r)},\
\forall t,r>0.  \neq
\item It holds that
\beq\label{thm22b} \pp_\beta\left(\frac 1t \int_0^t
g(X_s)ds>\alpha^{-1}(x)\right)\le \|\frac{d\beta}{d\mu}\|_2e^{-t
x},\ \forall t,x>0. \neq
\item For any $\lambda>0$,
\beq\label{thm22c}\Lambda(\lambda g):=\sup\left\{\int_\XX \lambda  g
h^2 d\mu -\EE(h,h)| h\in \dd(\EE), \mu(h^2)=1 \right\}\le
\alpha^*(\lambda)\neq where
$\alpha^*(\lambda):=\sup_{r\ge0}\{\lambda r - \alpha(r)\}$ is the
(semi)-Legendre transformation of $\alpha$.
\end{enumerate}
\nthm It is not completely contained in \cite[Theorem 2.4]{GLWY}
(the condition (A2) therein is not satisfied), but the proof there
works. Indeed $(a)\Leftrightarrow(b)$ and $(c)\Leftrightarrow(d)$
are obvious. We give the proof of the crucial implication
$(a)\implies (c)$ for its simplicity. In fact by the
transportation-information inequality in (a), we have for $r>0$,
$$
I(r)=\inf\{I(\nu|\mu);\ \nu(|g|)<+\infty,\ \nu(g)=r\}\ge \alpha(r)
$$
and then $I(r-)\ge \alpha(r)$ by the left-continuity of $\alpha$.
Hence the concentration inequality (\ref{thm22a}) follows
immediately from (\ref{prop21a}).

\brmk\label{rem21}{\rm By Rayleigh's principle, $\Lambda(\lambda g)$
is the supremum of the spectrum of the Schr\"odinger operator $\LL +
\lambda g$ (in the sum-form sense). }\nrmk

Bernstein's inequality (\ref{Bern1}) is just (\ref{thm22a}) with
$$
\alpha(r)= 1_{r\ge0}
\frac{2r^2}{\sigma^2\left(\sqrt{1+\frac{2Mr}{\sigma^2}}+1\right)^2}.
$$
Since $\alpha^{-1}(x)=\sqrt{2\sigma^2 x} + Mx$ for $x\ge0$, by
Theorem \ref{thm22}, Bernstein's inequality (\ref{Bern1}) is
equivalent  to  \beq\label{key1} \nu(g)\le \sqrt{2\sigma^2 I} + MI,\
I:=I(\nu|\mu), \ \forall \nu\in \MM_1(\XX) \text{ so that }\
\nu(|g|)<+\infty.\neq That is the strategy of this work.

Now let us present a very simple proof of Lezaud's result, which
illustrates also the main idea for our approaches to establish
(\ref{key1}). Assume $g\in L^2_0(\mu)$ so that $g^+\in
L^\infty(\mu)$.

Let $\nu=f\mu$ and $h=\sqrt{f}\in \dd(\EE)$ (trivial otherwise for
$I=+\infty$) such that $\nu(|g|)<+\infty$. Our main idea resides in
the following simple but key decomposition : \beq
\label{key2}\aligned \nu(g)&=\int_\XX g h^2 d\mu = \int_\XX g
\left[(h-\mu(h))^2 + 2 \mu(h) h\right] d\mu\ \ \ (\text{since }\mu(g)=0)\\
&=2\mu(h)\<g, h\>_\mu + \int_\XX g(h-\mu(h))^2 d\mu=: A + B.
\endaligned
\neq
{\it Bounding $A$.}\\
 For the first term $A=2\mu(h)\<g, h\>_\mu$, note that
$\mu(h)\le \sqrt{\mu(h^2)}=1$. Let $(-\LL)^{-1}g=\int_0^{+\infty}P_t
g dt$ be the Poisson operator (the integral is absolutely convergent
in $L^2(\mu)$ for all $g\in L^2_0(\mu)$ by the Poincar\'e
inequality). Hence
$$
\sigma^2 =\sigma^2(g)= 2 \int_0^\infty \<P_t g, g\> dt = 2
\<(-\LL)^{-1}g,g\>_\mu.
$$
By Cauchy-Schwarz, we have \beq\label{A} |\<g, h\>_\mu| \le
\sqrt{\<(-\LL)^{-1}g,g\> \EE(h,h)}= \sqrt{\frac {\sigma^2}{2} I}
\neq Hence $|A|\le \sqrt{2\sigma^2 I}$, in other words, the term $A$
is always bounded by the first term at the right hand side of the
inequality (\ref{key1}).

\brmk\label{rem22}  Even without the hypothesis of the Poincar\'e
inequality, (\ref{A}) is still true for $g\in L^2_0(\mu)$ by
Kipnis-Varadhan \cite{KV86} once if $\sigma^2(g)=2\int_0^\infty \<g,
P_tg\> dt<+\infty$. The latter condition is the famous sufficient
condition of Kipnis-Varadhan for the CLT of $\int_0^t g(X_s)ds$.
\nrmk

\noindent{\it Bounding $B$}.\\
Now for (\ref{key1}) it remains to prove that the second term $B$
satisfies \beq\label{key3}B=\int_\XX g[h-\mu(h)]^2 d\mu\le M
\EE(h,h)=MI. \neq It is indeed very easy in terms of $\|g\|_\infty$
: letting $g^+=\max\{g, 0\}$, we have by Poincar\'e,
$$
B=\int_\XX g[h-\mu(h)]^2 d\mu\le \int_\XX g^+[h-\mu(h)]^2 d\mu\le
\|g^+\|_\infty \Var_\mu(h)\le c_P \|g^+\|_\infty I.
$$
In other words we have proven (\ref{key1}) with $M=c_P
\|g^+\|_\infty$, which is a little better than Lezaud's estimate
$M=c_P \|g\|_\infty$. We summarize the discussion above as

\bprop\label{prop23} Let $g\in b\BB$ with $\mu(g)=0$. Then
 (\ref{key1}) holds with $M=c_P
\|g^+\|_\infty$, or equivalently Bernstein's inequality
(\ref{Bern1}) holds with such $M$. \nprop

Our remained task consists in proving (\ref{key3}) with some
constant $M=M(g)$ for various classes of functions $g$ under
different ergodicity conditions for the process. Remark that the
best constant $M(g)$ for (\ref{key3}) (or (\ref{key1})) is
positively homogeneous, i.e. $M(cg)=cM(g)$ for all $c\ge 0$.


\subsection{Approach by transportation-information inequality $T_cI$}

Let us introduce our first approach by means of the
transportation-information inequality $T_cI$ in \cite{GLWY}.

Consider a cost function $c: \XX^2\to [0,+\infty]$ which is always
lower semi-continuous (l.s.c.) and $c(x,x)=0$ for all $x\in\XX$,
here $c(x,y)$ represents the cost of transporting a unit mass from
$x$ to $y$. Now given two probability measures $\nu,\mu\in
\MM_1(\XX)$, we define the {\it transportation cost from $\nu$ to
$\mu$} by

\beq\label{cost1} T_c(\nu,\mu) :=\inf_{\pi\in \CC(\nu,\mu)}
\iint_{\XX^2} c(x,y) \pi(dx,dy) \neq where $\CC(\nu,\mu)$ is the
family of all couplings of $(\nu,\mu)$, i.e. all probability
measures $\pi$ on $\XX^2$ such that $\pi(A\times \XX)=\nu(A), \
\pi(\XX\times B)=\mu(B)$ for all $A,B\in\BB$.

Let $d(x,y)$ be a l.s.c. metric on $\XX$, which does not necessarily
generate the topology of $\XX$. For any $p\ge 1$, the quantity
\beq\label{Wasserstein1}
W_{p,d}(\nu,\mu):=\left(T_{d^p}(\nu,\mu)\right)^{1/p}=
\left(\inf_{\pi\in \CC(\nu,\mu)} \iint_{\XX^2} d^p(x,y) \pi(dx,dy)
\right)^{1/p}\neq is the so called {\it $L^p$-Wasserstein distance}
between $\nu$ and $\mu$. $W_{p,d}$ is a metric on
$\MM_{1}^{d,p}(\XX):= \{\nu\in \MM_1(\XX);~\left(\int_\XX
d^p(x_0,x)\nu(dx)\right)^{1/p}<+\infty\}$ ($x_0\in \XX$ is some
fixed point). We refer to the recent books of Villani \cite{Vill03,p-Vill05} for more on this subject.

An important particular case is $d(x,y)=1_{x\ne y}$, the trivial
metric on $\XX$. In that case

\beq\label{trivial} W_{1,d}(\nu,\mu)=\frac 12
\|\nu-\mu\|_{TV}=\sup_{A\in \BB}|\nu(A)-\mu(A)|\neq where
$\|m\|_{TV}=\sup_{f\in b\BB, |f|\le 1}|m(f)|$ is the total variation
of a signed bounded measure $m$ on $\XX$. More generally given a
positive continuous weight function $\phi$, consider the distance
$d_\phi(x,y)=1_{x\ne y} [\phi(x)+\phi(y)]$, then (cf. \cite{GL})
$$
W_{1,d_\phi}(\nu,\mu)=\|\phi(\nu-\mu)\|_{TV}.
$$
\bthm\label{thm23} Assume the following transportation-information
inequality \beq\label{thm23a} \alpha(T_c(\nu,\mu))\le I(\nu|\mu), \
\forall \nu\in \MM_1(\XX)\neq where $\alpha$ is nonnegative,
nondecreasing convex and left continuous with $\alpha(0)=0$ such
that its right inverse $\alpha^{-1}$ is concave and
$\alpha^{-1}(0)=0$. Then for every measurable $g\in L^2_0(\mu)$ such
that its sup-convolution

\beq\label{thm23b} g^*(y)=\sup_{x\in\XX} \left(g(x)-c(x,y)\right),\
y\in\XX \neq is in $L^1(\mu)$,  (\ref{key1}) and Bernstein's
inequality (\ref{Bern1}) hold with \beq\label{thm23c} M(g)=\mu(g^*)
c_P + c_P \alpha^{-1}\left(\frac 1{c_P}\right). \neq
 In particular if the $W_1I$-transportation-information inequality below holds
\beq\label{W1I} W_{1,d}^2(\nu,\mu)\le 2c_G I(\nu|\mu),\ \forall
\nu\in \MM_1(\XX)\neq then (\ref{key1}) holds for every
$d$-Lipschitzian function $g$ (with $\mu(g)=0$) with
$$
M(g)=\|g\|_{Lip(d)}\sqrt{2c_Pc_G}.
$$
\nthm

\bprf At first $g^*(y)\ge g(y), \ y\in \XX$, so $\mu(g^*)\ge
\mu(g)=0$. For (\ref{key1}) we may assume that $\nu=h^2 \mu$ with
$0\le h\in\dd(\EE)$ and $Var_\mu(h)\ne 0$ (trivial otherwise for
$\nu=\mu$). Letting $\tilde h=h-\mu(h)$ and $\tilde \nu:=\tilde
h^2\mu/\Var_\mu(h)$, we have by the very definition of $T_c$,
$$\aligned\int_\XX g(x) \tilde \nu(dx)&\le \int_\XX g^*(y) \mu(dy) +
T_c(\tilde\nu,\mu)\\
&\le \mu(g^*) + \alpha^{-1}(I(\tilde \nu|\mu))\le \mu(g^*) +
\alpha^{-1}\left(\frac{\EE(h,h)}{\Var_\mu(h)}\right)
\endaligned$$
where we have used $\EE(|\tilde h|, |\tilde h|)\le \EE(\tilde
h,\tilde h)=\EE(h,h)$. It follows by the concavity of $\alpha^{-1}$,
$$\aligned
B&=\int_\XX g\tilde h^2 d\mu\le \mu(g^*) \Var_\mu(h) +
\Var_\mu(h)\alpha^{-1}\left(\frac{\EE(h,h)}{\Var_\mu(h)}\right)\le
\mu(g^*) c_P I + c_P I \alpha^{-1}(1/c_P)
\endaligned$$
the desired (\ref{key3}).

For the last particular case we may assume that $\|g\|_{Lip(d)}=1$.
In that case $g^*=g$, and then one can apply (\ref{thm23c}).
 \nprf

\brmk{\rm By the preceding result, one can apply the criteria for
$T_cI$ or $W_1I$-transportation information inequalities in
\cite{GLWY} to obtain Bernstein's inequality.  }\nrmk


\section{Functional inequalities approach}

\subsection{Log-Sobolev inequality}
Recall that for $0\le f\in L^1(\mu)$, the entropy of $f$ w.r.t.
$\mu$ is defined by \beq\label{entropy} \Ent_\mu(f)=\mu(f\log
f)-\mu(f) \log \mu(f).\neq The log-Sobolev inequality (\cite{Ba92,
Led01}) says

\beq\label{logS} \Ent_\mu(h^2) \le 2c_{LS} \EE(h,h),\ \forall h\in
\dd(\EE), \neq where $c_{LS}$ is the best constant, called {\it
log-Sobolev constant. } It is well known that $c_P\le c_{LS}$.

\bthm\label{thm31} Assume the log-Sobolev inequality (\ref{logS}).
Let $g\in L^2_0(\mu)$ satisfy $\dsp \Lambda(\lambda):=\log \int_\XX
e^{\lambda g}d\mu<+\infty$ for some $\lambda>0$.

Then the transportation-information inequality (\ref{key1}) holds
with \beq\label{thm31a} M= \inf_{\lambda>0} \frac 1\lambda \left[c_P
\Lambda(\lambda) + 2 c_{LS}\right]\le c_P (\Lambda^*)^{-1}(\frac
{2c_{LS}}{c_P})\neq where $\Lambda^* : \rr^+\to [0,+\infty]$ is the
Legendre transform of $\Lambda$ and $(\Lambda^*)^{-1}$ is the right
inverse. In particular Bernstein's inequality (\ref{Bern1}) holds
with this constant $M$. \nthm

\bprf We may assume that  $\nu=h^2\mu$ with $0\le h\in \dd(\EE)$. We
have to bound the term $B=\int_\XX g[h-\mu(h)]^2 d\mu$ in the
decomposition (\ref{key2}). Writing $\tilde h=h-\mu(h),
I=I(\nu|\mu)=\EE(h,h)$, we have for any constant $\lambda>0$ such
that $\Lambda(\lambda)<+\infty$, $\int e^{\lambda g -a} d\mu=1$
where $a=\Lambda(\lambda)\ge0$, and then
$$
\aligned B &= \frac 1\lambda\left(\int_\XX (\lambda g - a) \tilde
h^2
d\mu + a \int \tilde h ^2 d\mu \right)\\
&\le \frac 1\lambda\left(\Ent_\mu(\tilde h^2)
 + a c_P I  \right)\\
&\le \frac 1\lambda\left[2c_{LS}
 + \Lambda(\lambda) c_P \right]\cdot I\\
\endaligned
$$
where the second inequality relies on $\Ent_\mu(f)=\sup_{g:
\mu(e^g)\le 1} \int_\XX fg d\mu$ (Donsker-Varadhan's variational
formula) and the Poincar\'e inequality, and the third one on the
log-Sobolev inequality. Optimizing over $\lambda>0$ yields
(\ref{key1}) with $M$ given in (\ref{thm31a}). \nprf

It is a surprise : the explicit estimate of $M=M(g)$ above is not
available even in the i.i.d. case under the exponential
integrability condition.

Let us give a more explicit estimate of $M$ in the diffusion case.
We assume that
\vskip5pt\noindent $(H_\Gamma)$ $(\EE, \dd(\EE))$ is
given by the carr\'e-du-champs $\Gamma: \dd(\EE)\times \dd(\EE)\to
L^1(\mu)$ (symmetric, bilinear definite nonnegative form):
\begin{equation}\label{Gamma}
\EE(h,h)=\int_\XX \Gamma(h,h)\,d\mu,\ \forall h\in\dd(\EE).
\end{equation}

\vskip5pt\noindent {\bf Diffusion framework. } We shall assume that
$\Gamma$ is a differentiation (or equivalently the sample paths of
$(X_t)$ are continuous, $\pp_\mu-a.s.$, cf. Bakry \cite{Ba92}), that
is: for all $(h_k)_{1\le k\le n}\subset \dd(\EE), g\in \dd(\EE)$ and
$F\in C_b^1(\rr^n)$,
$$
\Gamma(F(h_1,\cdots, h_n), g) = \sum_{i=1}^n
\partial_iF(h_1,\cdots, h_n)\Gamma(h_i, g).
$$
Write $\Gamma(f):=\Gamma(f,f)$ simply.

\bcor\label{cor31} Assume $(H_\Gamma)$ and that $\Gamma$ is a
differentiation. If the log-Sobolev inequality holds, then for any
$g\in \dd(\EE)$ so that $\Gamma(g)$ is bounded and $\mu(g)=0$, the
transportation-information inequality (\ref{key1}) holds with
\beq\label{thm31acor} M= 2c_{LS}\sqrt{c_P \|\Gamma(g)\|_\infty}.\neq

\ncor \bprf By Ledoux \cite{Led01} or Bobkov-G\"otze \cite{BG99}, in
the actual diffusion case the log-Sobolev inequality implies that
$$
\Lambda(\lambda)=\log \int_\XX e^{\lambda g}d\mu \le \frac 12
c_{LS}\lambda^2 \|\Gamma(g)\|_\infty, \ \forall \lambda>0.
$$
Plugging it into (\ref{thm31a}), we get $M\le 2c_{LS}\sqrt{c_P
\|\Gamma(g)\|_\infty}$. \nprf

\bexa\label{exa-OU}{\bf (Ornstein-Uhlenbeck processes) }{\rm  Let
$\mu=\NN(0,\theta)$, the Gaussian measure with zero mean and
variance $\theta>0$ on $\XX=\rr$, and $\LL f=f^{\prime\prime} -
\theta^{-1} x\cdot f'$. It is well known that $c_P=c_{LS}=\theta$.

For every Lipschitzian function $g$ with $\mu(g)=0$,
$\sqrt{\|\Gamma(g)\|_\infty}=\|\nabla g\|_\infty=\|g\|_{Lip}$ (the
Lipschitzian coefficient w.r.t. the Euclidean metric). By Corollary
\ref{cor31}, Bernstein's inequality (\ref{Bern1}) holds with
$M=2c_{LS}\sqrt{c_P}\|g\|_{Lip}=2\theta^{3/2}\|g\|_{Lip}$. It is
worth mentioning that for the special observable $g(x)=x$,
(\ref{key1}) and then Bernstein's inequality (\ref{Bern1}) hold with
$M=0$ (i.e. the corresponding Gaussian concentration inequality
holds); and for general $g$ with $\mu(g)=0$,
$$
\nu(g) \le \|g\|_{Lip}\sqrt{2 \theta I}
$$
holds by \cite[Proposition 2.9]{GLWY}.

 But by Theorem \ref{thm31}, for every $\mu$-centered function
$g$ such that $\int e^{\delta g} d\mu<+\infty$ (for instance if
$g\le C (1+|x|^2)$), Bernstein's inequality (\ref{Bern1}) holds with
$M=M(g)$ given in (\ref{thm31a}). Though natural, that was not known
before up to our knowledge. It is easy to see that Bernstein
inequality is false for observable $g(x)$ such that
$\lim_{x\to\infty} \frac {g(x)}{x^2}=+\infty$.

Let us look at the particularly interesting observable
$g(x)=g_0(x):=x^2 -\theta$ for which we can get sharp Bernstein
inequality. Indeed since $-\LL g_0=-2\theta^{-1} g_0$,
$$\sigma^2(g_0)= 2\<(-\LL)^{-1}g_0,g_0\>_\mu= \theta
\Var_\mu(g_0)=2 \theta^3.$$ On the other hand observe that for each
real number $a<\frac 1{2}$, $U(x):=\exp\left(\frac{ax^2}{2
\theta}\right)\in L^2(\mu)$, and
$$
\left[\LL  + \frac{a-a^2}{\theta^{2}} g_0\right] U=
\frac{a^2}{\theta} U.
$$
In other words $U$ is a positive eigenfunction of the Schr\"odinger
operator $\LL  + \frac{a-a^2}{\theta^2}g_0$ associated with
eigenvalue $a^2/\theta$, which implies that (by Perron-Frobenius
theorem and Rayleigh's formula)
$$
\Lambda\left(\frac{a-a^2}{\theta^2}g_0\right)=\frac{a^2}{\theta},\
a<\frac 1{2}.
$$
Hence for all $\lambda<\lambda_0:=\frac 1{4\theta^2}$, taking $a=a_-
:= \frac 12\left(1- \sqrt{1-4\theta^2\lambda}\right)<1/2$, we have
$$
\Lambda(\lambda g_0)= \frac 1{4\theta}\left(1-
\sqrt{1-4\theta^2\lambda}\right)^2
$$
Since $\lambda\to \Lambda(\lambda g_0)$ from $\rr$ to
$(-\infty,+\infty]$ is convex and lower semi-continuous, and its
left derivative at $\lambda_0$ is $+\infty$, we conclude that
\beq\label{OU1} \Lambda(\lambda):=\Lambda(\lambda g_0)=\frac
1{4\theta}\left(1- \sqrt{1-4\theta^2\lambda}\right)^2\ \text{ if }\
\lambda\le \lambda_0=\frac 1{4\theta^2};\ +\infty, \ \text{ if
}\lambda>\lambda_0. \neq From the previous explicit expression we
obtain (by the fact that the geometric mean is not greater than the
arithmetic mean)
$$
\Lambda(\lambda)= \frac{\sigma^2(g_0)\lambda^2 }{2 (\frac
12[1+\sqrt{1-4\theta^2\lambda}]^2}\le \frac{\sigma^2(g_0)\lambda^2
}{2(1-4\theta^2\lambda)}, \ \lambda\in (0, \lambda_0)
$$
where it follows that $g_0(x)=x^2-\theta$ satisfies the Bernstein
inequality (\ref{Bern1}) with the sharp constant $M=4\theta^2$.

Notice that (\ref{OU1}) will give, by Theorem \ref{thm22}, the
concentration inequality for the estimator $\frac 1t\int_0^t X_s^2
ds$ of $\theta$, which is not only sharp for the CLT and moderate
deviation scales, but also for large deviations.
 }\nexa


\subsection{$\Phi$-Sobolev inequality}

Let $\Phi:\rr^+\to [0,+\infty]$ be a Young function, i.e. a convex,
increasing and left continuous function with $\Phi(0)=0$ and
$\lim_{x\to+\infty} \Phi(x)=+\infty$. Consider the Orlicz space
$L^{\Phi}(\mu)$ of those measurable functions $g$ on $\XX$ so that
its gauge norm

$$
N_{\Phi}(g):=\inf\{c>0; \int \Phi(|g|/c) d\mu\le 1\}
$$
is finite, where the convention $\inf \emptyset:=+\infty$ is used.
The Orlicz norm of $g$ is defined by
$$
\|g\|_\Phi:=\sup\{\int g u \,d\mu; \ N_{\Psi}(u)\le 1\}
$$
where \beq \Psi(r):=\sup_{\lambda\ge 0} (\lambda r -
\Phi(\lambda)),\ r\ge0 \neq is the convex conjugate of $\Phi$. It is
well known that (\cite[Proposition 4, p.61]{RR})
$$
N_\Phi(g)\le \|g\|_\Phi\le 2 N_\Phi(g).
$$

The $\Phi$-Sobolev inequality says that \beq \label{PhiS}
\|(h-\mu(h))^2\|_{\Phi} \le c_{P,\Phi} \EE(h,h), \ \forall
h\in\dd(\EE) \neq called sometimes Orlicz-Poincar\'e inequality,
where $c_{P,\Phi}$ is the best constant. There is a rich theory of
long history for this subject, see \cite{Chen05, Led01, Wang}.

Set $\tilde \Phi(x):=\Phi(x^2), x\ge0$ and let $\tilde \Psi$ be the
Legendre transform of $\tilde \Phi$.

\blem\label{lem34} Assume the $\Phi$-Sobolev inequality
(\ref{PhiS}). If $g\in L^{\tilde \Psi}(\mu)$ so that $\mu(g)=0$,
then $\int_0^t g(X_s) ds \in L^2(\pp_\mu)$ and it holds that
\beq\label{lem34aa} \sigma^2(g)=\lim_{t\to+\infty} \frac 1t {\rm
Var}_{\pp_\mu} \left(\int_0^t g(X_s) ds\right)\le c_{P,\Phi}
\|g\|^2_{\tilde \Psi}. \neq Moreover \beq\label{lem34b}
\<g,h\>^2_\mu \le \frac 12 \sigma^2(g) \EE(h,h), \ \forall h\in
\dd(\EE). \neq
 \nlem

\bprf At first for $g\in L^2_0(\mu)$, notice that by the spectral
decomposition and Cauchy-Schwarz,
$$
\sqrt{\<g, (-\LL)^{-1}g\>_\mu} =\sup_{h\in \dd(\EE), \EE(h,h)\le 1}
\<g,h\>_\mu
$$
and
$$
|\<g,h\>_\mu| =|\<g,h-\mu(h)\>_\mu|\le \|g\|_{\tilde \Psi} N_{\tilde
\Phi} (h-\mu(h)).
$$
Furthermore by the $\Phi$-Sobolev inequality (\ref{PhiS}),
$$N_{\tilde \Phi} (h-\mu(h))=
\sqrt{N_{\Phi}((h-\mu(h))^2)}\le \sqrt{ \|(h-\mu(h))^2\|_\Phi}\le
\sqrt{c_{P,\Phi}\EE(h,h)}$$ therefore \beq\label{lem34a}\<g,
(-\LL)^{-1}g\>_\mu \le c_{P,\Phi} \|g\|_{\tilde \Psi} ^2,\ g\in
L^2_0(\mu). \neq

Now take a sequence $(g_n)$ in $L^\infty_0(\mu)$ converging to $g$
in $L^{\Psi}(\mu)$, we have for any $t>0$,
$$\aligned
\frac 1t {\rm Var}_{\pp_\mu}\left(\int_0^t (g_n-g_m)(X_s)
ds\right)&\le \sigma^2(g_n-g_m)=2 \<g_n-g_m,
(-\LL)^{-1}(g_n-g_m)\>_\mu\\
&\le 2 c_{P,\Phi} \|g_n-g_m\|_{\tilde \Psi}^2. \endaligned$$ This
implies not only ``$\int_0^t g(X_s) ds \in L^2(\pp_\mu)$" but also
(\ref{lem34aa}). The last claim (\ref{lem34b}) holds for $g_n$ in
place of $g$ then remains true for $g$ by letting $n\to\infty$.
\nprf

\bthm \label{thm32} Assume the $\Phi$-Sobolev inequality
(\ref{PhiS}) and let $\Psi$ be the convex conjugate of $\Phi$ given
above. If $g\in L^{\tilde \Psi}(\mu)$ and $g^+\in L^{\Psi}(\mu)$
with $\mu(g)=0$, then the transportation-information inequality
(\ref{key1}) holds with $\sigma^2=\sigma^2(g)$ given by
(\ref{lem34aa}) and \beq\label{thm32a} M = N_\Psi(g^+) \cdot
c_{P,\Phi}.\neq In particular Bernstein's inequality (\ref{Bern1})
holds with that constant $M$. \nthm

\bprf The proof is even easier than that of Theorem \ref{thm31}. For
(\ref{key1}) we may assume that $\nu=h^2\mu$ with $0\le h\in
\dd(\EE)$. By Lemma \ref{lem34}, $\sigma^2=\sigma^2(g)$ given by
(\ref{lem34aa}) is finite.  The term $A$ in (\ref{key1}) is bounded
by $\sqrt{2\sigma^2 I}$ by (\ref{lem34b}). For the term $B=\int_\XX
g[h-\mu(h)]^2 d\mu$ we have
$$
B \le N_\Psi(g^+) \|[h-\mu(h)]^2\|_\Phi \le c_{P,\Phi}  N_\Psi(g^+)
I
$$
where the desired result follows. \nprf

\brmk{\rm When $\Phi(x)=|x|$, $\Psi(x)=+\infty\cdot 1_{x>1}$,
$N_\Psi(h)=\|h\|_\infty$. Then this result generalizes Proposition
\ref{prop23}. }\nrmk

\brmk{\rm For one-dimensional diffusions, an explicit necessary and
sufficient condition for the  $\Phi$-Sobolev inequality (\ref{PhiS})
is available, see the book of M.F. Chen \cite{Chen05}. For
$\Phi$-Sobolev inequality in high dimension, see the book of F.Y.
Wang \cite{Wang} for numerous known results.}\nrmk

\bexa\label{exa31} {\rm As a well known fact (see Saloff-Coste
\cite{LSC}), for the Brownian Motion $(B_t)$ on a compact connected
Riemannian manifold $M$ of dimension $n$ with the invariant measure
$\mu$ given by the normalized Riemannian measure $\dsp
\frac{dx}{V(M)}$(where $V(M)$ is the volume of $M$), the Dirichlet
form $\int |\nabla f|^2d\mu$ satisfies the $\Phi$-Sobolev inequality
(\ref{PhiS}) with
$$
\Phi(t)=\left\{\begin{array}{cc}
+\infty I_{(1,\infty)}(|t|), &{\rm if}  \quad{n=1, }\\
\hbox{}\\
\dsp \exp(C |t|) -1, &{\rm if}   \quad{n=2, }\\
\hbox{}\\
|t|^\frac{2n}{n-2},&{\rm if}  \quad{n\ge 3.}
\end{array}
\right.
$$
Hence Bernstein's inequality (\ref{Bern1}) holds for $g\in
L^1_0(\mu)$ satisfying
$$
g\in\left\{\begin{array}{cc}
L^1(\mu), &{\rm if}  \quad{n=1, }\\
\hbox{}\\
L^1\log L^1, &{\rm if}   \quad{n=2, }\\
\hbox{}\\
L^\frac{2n}{n+2}(\mu),&{\rm if}  \quad{n\geq 3.}
\end{array}
\right.
$$
Those still hold for diffusion generated by $\Delta -\nabla
V\cdot\nabla$ with $C^2$-smooth function  $V$ on a connected compact
manifold. } \nexa

\bexa\label{exa32} {\rm  Consider the measure $\dsp \mu_\beta(dx)=
\frac{\exp(-|x|^\beta)}{Z_\beta}$ (where $Z_\beta$ is the normalized
constant), and $\beta>1$. For the diffusion process corresponding to
the Dirichlet form $\<-\LL f, f\>_\mu=\int |\nabla f|^2d\mu$, it
satisfies $\Phi$-Sobolev inequality (\ref{PhiS}) with
$$\Phi_\alpha(x)=x \log^\alpha
(1+x),\ \alpha=2(1-1/\beta)$$according to Barthe, Cattiaux and
Roberto \cite[section 7]{BCR}. Hence Bernstein's inequality
(\ref{Bern1}) holds for $g\in L^2_0(\mu)$ satisfying
\beq\label{exa32-a} \int \exp\left(\lambda
(g^+)^{\beta/(2\beta-2)}\right) d\mu <+\infty, \ \text{for some }\
\lambda>0. \neq }\nexa

Those two examples show that for Bernstein's inequality to hold, the
integrability condition on the observable $g$ in the continuous time
symmetric Markov processes case may be much weaker than the
exponential integrability condition in the i.i.d. case.


\section{Lyapunov function method}
Sometimes functional inequalities are difficult to check. In that
situation the easy-to-check Lyapunov function method will be very
helpful.

\subsection{General result}
A measurable function $G$ is said to be in the $\mu$-extended domain
$\dd_{e,\mu}(\LL)$ of the generator of the Markov process $((X_t),
\pp_\mu)$ if there is some measurable function $g$ such that
$\int_0^t |g|(X_s)\,ds<+\infty, \pp_\mu$-a.s.\! and one
$\pp_\mu$-version of
$$
M_t(G):=G(X_t) - G(X_0)+\int_0^t g(X_s)ds
$$
is a local $\pp_\mu$-martingale. It is obvious that $g$ is uniquely
determined  up to $\mu$-equivalence. In such case one writes
$G\in\dd_{e,\mu}(\LL)$ and $-\LL G=g$. When the above properties
hold for $\pp_x$ instead of $\pp_\mu$ for {\bf every } $x\in\XX$, we
say that $G$ belongs to the extended domain $\dd_e(\LL)$. In the
latter case $-\LL G=g$ is determined uniquely up to $\int_0^\infty
e^{-t} P_t(x,\cdot) dt$-equivalence for every $x\in \XX$.

The Lyapunov condition can be stated  now :
\begin{itemize}
\item[$(H_L)$] There exist a measurable function $U:\mathcal{X}\to
[1,+\infty)$ in $\dd_{e,\mu}(\LL)$, a positive function $\phi$
and a constant $b>0$ such that
$$- \frac{\mathcal{L}U}{U}\ge \phi -b, \ \mu\textrm{-a.s.}$$
\end{itemize}
When the process is irreducible and the constant $b$ is replaced by
$b1_C$ for some ``small set" $C$, then it is well-known that the
existence of a positive  bounded $\phi$  such that
$\inf_{\XX\setminus C}\phi>0$ in $(H_L)$ is equivalent to Poincar\'e
inequality (see \cite{BBCG, BCG2}, for instance).

Lyapunov conditions are widely used to study the speed of
convergence of Markov chains \cite{MT} or Markov processes
\cite{DMT,DFG}, large or moderate deviations and essential spectral
radii \cite{Wu01, G01, Wu04} or sharp large deviations \cite{KM}. More
recently, they have been used to study functional inequalities such as
weak Poincar\'e inequality \cite{BCG2} or super-Poincar\'e
inequality \cite{CGWW}. See Wang
 \cite{Wang} on weak and super Poincar\'e inequalities.

For a given function $f$, let $K_\phi(f)\in [0,+\infty]$ be the
minimal constant $C\in [0,+\infty]$ such that $|f|\le C \phi$.

\bthm \label{thm51} Assume the Lyapunov function condition $(H_L)$.
For $g\in L^2_0(\mu)$, if $K_\phi(g^+)<+\infty$, then the
transportation-information inequality (\ref{key1}) holds with
\beq\label{thm32aly} M = K_\phi(g^+)\left(bc_P +1\right).\neq In
particular Bernstein's inequality (\ref{Bern1}) holds with that
constant $M$. \nthm

\bprf
 We are inspired by the elegant proof of Barthe-Bakry-Cattiaux-Guillin \cite{BBCG}
 for the Poincar\'e inequality. As before let $\nu=h^2\mu$
with $0\le h\in \dd(\EE)$. For the term $B=\int_\XX g[h-\mu(h)]^2
d\mu$ in (\ref{key1}) we have by $(H_L)$,
$$
B\le K_\phi(g^+) \int_\XX \phi [h-\mu(h)]^2 d\mu\le K_\phi(g^+)
\int_\XX \left(b- \frac {\LL U}{U} \right) [h-\mu(h)]^2 d\mu.
$$
By a result in large deviations \cite[Lemma 5.6]{GLWY}, we have
$$
\int_\XX -\frac {\LL U}{U}[h-\mu(h)]^2 d\mu\le \EE(h,h)=I.
$$
Hence applying the Poincar\'e inequality, we get
$$
B\le K_\phi(g^+)\left(bc_P +1\right) I
$$
the desired result. \nprf

\subsection{Particular case : diffusions on $\rr^d$} Let $\XX=\rr^d$, $x\cdot y$ and $|x|=\sqrt{x\cdot x}$ be
the Euclidean inner product and norm, respectively. Consider
$\LL=\Delta -\nabla V\cdot \nabla$ on $\rr^d$,  where $V$ is lower
bounded $C^2$-smooth such that $Z=\int_{\rr^d} e^{-V} dx$ is finite.
The corresponding semigroup $P_t$ is symmetric on $L^2(\mu)$ for
$\mu=\frac 1Z e^{-V} dx$. From Theorem \ref{thm51} we derive easily

\bcor\label{cor1} In the framework above, let $\gamma>0$ be some
fixed constant. If one of the following conditions
\begin{equation}\label{kustr}
\exists a<1, R, c>0, {\rm such~that~if~ }\ |x|>R,\qquad (1-a)|\nabla
V|^2-\Delta V\ge c~(1+|x|^\gamma)
\end{equation}
or
\begin{equation}\label{simpl}
\exists R,c>0, {\rm such~that~}\forall |x|>R, \qquad |x|^{\gamma/2}
\frac x{|x|}\cdot\nabla V(x)\ge c~(1+|x|^\gamma)
\end{equation}
is satisfied,  then the Lyapunov function condition $(H_L)$ is
satisfied with $\phi(x):=c(1+|x|^\gamma)$, and then for any
$\mu$-centered function $g$ such that $g(x)\le C (1+|x|^\gamma)$,
Bernstein's inequality (\ref{Bern1}) holds for some constant
$M=M(g)$ given by (\ref{thm32aly}). \ncor

\begin{proof} Under (\ref{kustr}), one takes $U=e^{aV}$; and under
(\ref{simpl}) one choose $U=e^{a|x|^{1+(\gamma/2)}}$ with small enough
$a>0$ (so that $c$ may be arbitrary). One sees that condition $(H_\Gamma)$ is satisfied in both
cases.
\end{proof}

\bexa\label{exa51}{\rm Let $V(x)=|x|^\beta$ ($\beta>0$ is fixed) for
$|x|>1$ in the framework above.

{\bf Case 1. $\beta\in (0,1)$.} In this case the Poincar\'e
inequality does not hold (cf. \cite{Led01}). And Bernstein's
inequality (\ref{Bern1}) does not hold for {\bf all } $g\in b\BB$
(with $\mu(g)=0$) as explained in the Introduction. Section 6 is
devoted to such examples.

{\bf Case 2. $\beta=1$.} For this exponential type's measure $\mu$,
the Poincar\'e inequality holds and one can apply Lezaud's result
for bounded $g$. We do not believe that the Bernstein's inequality
holds for unbounded $g$.

{\bf Case 3. $\beta> 1$.} Condition (\ref{simpl}) is satisfied with
$\gamma=2(\beta-1)$. Hence Bernstein's inequality (\ref{Bern1})
holds for $\mu$-centered $g$ such that $g\le C(1+|x|^{2(\beta-1)})$,
in concordance with condition (\ref{exa32-a}) in Example
\ref{exa32}. }\nexa

\subsection{Particular case : birth-death processes}
Let $\XX=\nn$ and
$$
\LL f(k)=b_k (f(k+1)-f(k)) + a_k (f(k-1)-f(k)), \ k\in\nn
$$
where $b_k>0, k\ge 0$  are the birth rates, $a_k>0, k\ge 1$ are the
death rates respectively, and $f(-1):=f(0)$.

We assume that the process is positive recurrent, i.e.,
$$\sum_{n\ge 0}\pi_n\sum_{i\ge n}(\pi_i b_i)^{-1}=\infty \quad
\text{and} \quad C:=\sum_{n=0}^{+\infty}\pi_{n}<+\infty,$$ where
$\pi_n$ is given by
$$\pi_{0}=1, \quad \pi_{n}=\frac{b_{0}b_{1}\cdots
b_{n-1}}{a_{1}a_{2}\cdots a_{n}}, \quad n\geq 1$$ is an invariant
measure of the process. Define the normalized probability $\mu$ of
$\pi$ by $\mu_{n}=\frac{\pi_{n}}{C}$ for any $n\geq0,$ which is
actually the unique reversible invariant probability of the process.

\bcor\label{cor52} Given a positive weight function $\phi_0$ on
$\nn$ such that $\phi_0\ge \delta>0$. If there are some constant
$\kappa>1$ and some $N\ge 1$ so that \beq\label{cor52a} a_n - \kappa
b_n\ge \phi_0(n), \ n\ge N, \neq then $(H_\Gamma)$ holds with
$\phi(n):=(1-\kappa^{-1})\phi_0(n)$ (and some finite constant $b$).
In particular the results in Theorem \ref{thm51} holds true. \ncor

\bprf Let $U(n)=\kappa^n$, we have
$$
-\frac{\LL U}{U} (n)=\frac{\kappa-1}{\kappa}\left(a_n - \kappa b_n
\right)
$$
where it follows that $c_P<+\infty$ (\cite{BBCG, BCG2}) and so the
desired result holds by Theorem \ref{thm51}. \nprf

\bexa\label{exa52}{\bf ($M/M/\infty$-queue system) }{\rm Let
$b_k=\lambda>0$ ($k\ge 0$) and $a_k=k$ ($k\ge 1$). Then $\mu$ is the
Poisson distribution with parameter $\lambda$. It is an ideal model
for a queue system with a number of {\it serveurs} much larger than
the number of clients. It is well known that $c_P=1$ but the
log-Sobolev inequality does not hold (\cite{Wu00h}).

For $\phi_0(n)=n+\delta$ where $\delta>0$ is fixed, taking
$U(n)=\kappa^n$ ($\kappa>1$) as above and applying Theorem
\ref{thm51}, we get by an optimization over $\kappa>1$ that for all
$g$ so that $g\le K(n+\delta)$ ($K>0$), $B\le MI$ where
\beq\label{exa52a} M=K[(\sqrt{\lambda}+1)^2 +  \delta]. \neq Hence
(\ref{key1}) and Bernstein's inequality (\ref{Bern1}) hold with such
$M$. Notice that the growth of $M$ for large $\lambda$ is linear in
$\lambda$.

An important observable is $g_0(n)=n-\lambda$ (then $L_t(g_0)$ is
the difference between the mean number of clients in the queue
system during time interval $[0,t]$ and the asymptotic mean
$\lambda$).  Since $(-\LL)^{-1} g_0=g_0$, we have
$\sigma^2(g_0)=2\<(-\LL)^{-1}
g_0,g_0\>_\mu=2\Var_\mu(g_0)=2\lambda$. We want to get a better
estimate of $M=M(g_0)$.

For $U(n)=\kappa^n$ $(\kappa>0)$, we have
$$\left[\LL + \frac{\kappa-1}{\kappa} g_0\right]U =\frac{(\kappa-1)^2}{\kappa} \lambda U.
$$
In other words $0<U\in L^2(\mu)$ is an eigenfunction of the
Schr\"odinger operator $\LL + \frac{\kappa-1}{\kappa} g_0$ with
eigenvalue $\frac{(\kappa-1)^2}{\kappa} \lambda$. By
Perron-Frobenius theorem and Raylaigh's principle,
$$
\Lambda\left(\frac{\kappa-1}{\kappa} g_0 \right)=
\frac{(\kappa-1)^2}{\kappa} \lambda.
$$
Thus if $s<1$, \beq\label{sharp} \Lambda(s g_0)=\frac{\lambda
s^2}{1-s}= \frac{\sigma^2(g_0) s^2}{2(1-s)} \neq and then
$\Lambda(sg_0)=+\infty$ for all $s\ge 1$ (by the convexity of $s\to
\Lambda(sg_0)$).

By Theorem \ref{thm22}, for $g=g_0$, not only the Bernstein
inequality (\ref{Bern1}) holds with the optimal constant $M(g_0)=1$,
and this inequality is itself sharp : indeed (\ref{sharp}) implies
by Proposition \ref{prop21} and the large deviation lower bound in
Wu \cite[Theorem B.1]{Wu00b},
$$
\lim_{t\to\infty} \frac 1t \log \pp_\mu\left( \frac 1t \int_0^t X_s
ds > \lambda + r \right) =
-\frac{r^2}{\lambda\left(\sqrt{1+\frac{r}{\lambda}}+1\right)^2},\
r>0.
$$
The calculus above shows that the mean number of clients $\frac 1t
\int_0^t X_s ds$ does not possess any Poisson type's concentration
inequality, contrary to the intuition that one might have for this
standard process related with the Poisson measure.
 } \nexa


\section{A Lipschitzian approach}
In this section we assume always the existence of
 the carr\'e-du-champs operator $\Gamma$, i.e. $(H_\Gamma)$ in \S3.
 We suppose furthermore that $\Gamma=\Gamma_0+\Gamma_1$ where
 $\Gamma_k: \dd(\EE^2)\to L^1(\mu),\ k=0,1$ are both
 bilinear nonnegative definite forms, $\Gamma_0$ is a
 differentiation, $\Gamma_1$ is given by
 $$\aligned
&\Gamma_1(f,g)(x)=\frac 12\int_\XX (f(y)-f(x))(g(y)-g(x)) J(x,dy),\
f,g\in \dd(\EE).
\endaligned$$
Here $\Gamma_0$ corresponds to the continuous diffusion part of
$(X_t)$, and  $J(x,dy)$ is a nonnegative jumps kernel (maybe
$\sigma$-infinite) on $\XX$ such that $J(x,\{x\})=0$ and
$\mu(dx)J(x,dy)$ is symmetric on $\XX^2$, describing the jumps rate
of the process.

\subsection{General result}
Recall that $\Gamma(f)=\Gamma(f,f)$.
 \bthm\label{thm41} Assume that $d$ is a lower semi-continuous
metric on $\XX$ (which does not necessarily generate the topology of
$\XX$), such that $\int_\XX d(x,x_0)^2 d\mu(x)<+\infty$. Given $g\in
L^2_0(\mu)$, let $G\in L^2_0(\mu)\bigcap \dd_2(\LL)$ be the unique
solution of the Poisson equation $-\LL G=g$. If
$\|\Gamma(G)\|_\infty<+\infty$, then the transportation-information
inequality (\ref{key1}) holds with \beq\label{thm41a} M = 2
\sqrt{c_P\|\Gamma(G)\|_\infty}.\neq In particular Bernstein's
inequality (\ref{Bern1}) holds with that constant $M$. \nthm

\bprf As before we may assume that $\nu=h^2\mu$ with $0\le h\in
\dd(\EE)\bigcap L^\infty(\mu)$. For the term $B=\int_\XX
g[h-\mu(h)]^2 d\mu$ in (\ref{key2}), setting $\tilde h=h-\mu(h)$ we
 write
$$\aligned
B &= \<-\LL G, \tilde h^2\>_\mu = \int_\XX \Gamma_0(G, \tilde
h^2)d\mu + \int_\XX \Gamma_1(G, \tilde h^2)d\mu.
\endaligned
$$
For the $\Gamma_0$-term, we have
$$\aligned
\int_\XX \Gamma_0(G, \tilde h^2)d\mu&\le \int_\XX
\sqrt{\Gamma_0(G)\Gamma_0( \tilde h^2)}d\mu=2 \int_\XX
\sqrt{\Gamma_0(G)\tilde h^2\Gamma_0(h)}d\mu
\endaligned$$
The $\Gamma_1$-term above requires some more work. We proceed as
follows.
$$\aligned
\int_\XX \Gamma_1(G, \tilde h^2)d\mu &= \frac 12
\iint_{\XX^2}(G(y)-G(x)) (\tilde h(y)+\tilde h(x)) (\tilde
h(y)-\tilde h(x)) \mu(dx)J(x,dy)\\
 &\le 2 \int_{\XX} \mu(dx) \sqrt{\int_\XX
(\tilde h(y)-\tilde h(x))^2 \mu(dx)J(x,dy)}\\
&\quad\quad\cdot \sqrt{\frac 18\int_{\XX}(G(y)-G(x))^2 [\tilde
h(y)+\tilde h(x)]^2 \mu(dx)J(x,dy)}.
\endaligned
$$
Plugging those two estimates into the expression of $B$ above, we
get by Cauchy-Schwarz's inequality,
$$\aligned
B&\le 2 \sqrt{\int_\XX \Gamma_0(G)\tilde h^2d\mu + \frac 18\int_\XX
\int_{\XX}(G(y)-G(x))^2 [\tilde h(y)+\tilde h(x)]^2
\mu(dx)J(x,dy)}\\
&\quad\quad\cdot \sqrt{\int_\XX (\Gamma_0(h)+\Gamma_1(h)) d\mu}.
\endaligned
$$
The last factor is $\sqrt{I}$. Using the symmetry in $(x,y)$ of
$\mu(dx)J(x,dy)$ and $(a+b)^2\le 2 (a^2 +b^2)$ , the second term
inside the first square root above can be bounded by
$$\aligned
\frac 14&\iint_{\XX^2}(G(y)-G(x))^2 [\tilde h(y)^2+\tilde h(x)^2]
\mu(dx)J(x,dy)\\
 &= \frac 12 \iint_{\XX^2}(G(y)-G(x))^2 \tilde h(x)^2
\mu(dx)J(x,dy)=  \int_\XX \Gamma_1(G)(x) \tilde h(x)^2 \mu(dx).
\endaligned$$
Hence the sum inside the first square root above is not greater than
$ \int_\XX \Gamma(G)(x) \tilde h(x)^2 \mu(dx)$. Thus we obtain
\beq\label{thm41b} B=\int_\XX g [h-\mu(h)]^2d\mu\le 2\sqrt{\int_\XX
\Gamma(G)(x) \tilde h(x)^2 \mu(dx)}\cdot \sqrt{I}.  \neq Now noting
that $ \int_\XX \Gamma(G)(x) \tilde h(x)^2 \mu(dx)\le
\|\Gamma(G)\|_\infty \Var_\mu(h)\le c_P \|\Gamma(G)\|_\infty I$, we
conclude that $B\le 2\sqrt{ c_P \|\Gamma(G)\|_\infty} I$, the
desired result. \nprf

Some sharp estimates of $\|\Gamma(G)\|_\infty$ for diffusions are
available : see Djellout and Wu \cite{DW} for one dimensional
diffusions, and Wu \cite{Wu09} for elliptic diffusions on manifolds.
Here we present examples of jumps processes.


\subsection{Birth-death processes continued}

The following two lemmas are taken from Liu and Ma \cite{LiuMa06}.

\blem\label{solu} Given a function $g$ on $\nn$ with $\mu(g)=0,$
consider the Poisson equation \beq\label{eqofQ} -\mathcal{L}G=g.
\neq For any $k\ge 0,$ the solution of the above equation
(\ref{eqofQ})
 satisfies the following relation :
 \begin{equation}\label{relation}
 G(k+1)-G(k)=-\frac{\sum_{j=0}^{k}\mu_jg(j)}{\mu_{k+1}a_{k+1}}=\frac{\sum_{j\ge k+1}\mu_jg(j)}{\mu_{k+1}a_{k+1}}.
 \neq\nlem

\blem\label{compari0} Let $\rho:\nn\to \rr$ be an increasing
function in $L^2(\mu)$. Provided that $||g||_{Lip(\rho)}:=\sup_{k\in
\nn}\frac{|g(k+1)-g(k)|}{\rho(k+1)-\rho(k)}=1$ with $\mu(g)=0$, we
have for any $k\ge 0$, \beq\label{compari1} \sum_{i\ge
k}\mu_ig(i)\le\sum_{i\ge k}\mu_i(\rho(i)-\mu(\rho)).\neq \nlem

We can derive easily

\bcor\label{cor42} Let $\rho:\nn\to \rr$ be an increasing function
in $L^2(\mu)$. If \beq\label{cor42a} K:=\frac 12\sup_{n\ge 0}
\left(1_{n\ge 1}\frac 1{a_n\mu_n^2} \left[\sum_{i\ge
n}\mu_i(\rho(i)-\mu(\rho))\right]^2+\frac{1}{b_n\mu_{n}^2}\left[\sum_{i\ge
n+1}\mu_i(\rho(i)-\mu(\rho))\right]^2 \right) \neq is finite, then
for every $g$ with $\mu(g)=0$ and $\|g\|_{Lip(\rho)}<+\infty$, the
transportation inequality (\ref{key1}) holds with $M=2\sqrt{c_P K}
\|g\|_{Lip(\rho)}$. \ncor

\bprf By Lemmas \ref{solu} and \ref{compari0}, the solution $G$ of
$-\LL G=g$ satisfies $\|\Gamma(G)\|_\infty\le K\|g\|_{Lip(\rho)}^2$
(using $a_{n+1}\mu_{n+1}=b_n\mu_n$). It remains to apply Theorem
\ref{thm41}. \nprf

See \cite{LiuMa06} for convex concentration inequalities. Though we
can give many examples to which Corollary \ref{cor42} applies, we
want to look at the $M/M/\infty$ queue system again.

\bexa{\bf ($M/M/\infty$ queue, continued) } {\rm The constant $K$ in
(\ref{cor42a}) above is infinite for $\rho(n)=n$, but finite for
$\rho(n)=\sum_{k=0}^n 1/\sqrt{k+1}$ (a quite artificial choice).
{\it What happens for $\rho(n)=\rho_0(n):=n$ ?} (In that case
$\|g\|_{Lip(\rho_0)}=:\|g\|_{Lip}$ is the Lipschitzian coefficient
w.r.t. the Euclidean metric.)

A crucial feature of this model is the commutation relation
$DP_t=e^{-t}P_tD$ where $Df(n):=f(n+1)-f(n)$, a property shared by
Ornstein-Uhlenbeck process for $D=\nabla$.  From this fact one sees
that
$$
\|(-\LL)^{-1}g\|_{Lip} \le \|g\|_{Lip}.
$$
Then if $\|g\|_{Lip}\le 1$, $G=(-\LL)^{-1}g$ satisfies
$$
\Gamma(G)(n)= \frac 12 \left(\lambda [G(n+1)-G(n)]^2 + n
[G(n-1)-G(n)]^2 \right)\le \frac 12 (\lambda+n).
$$
Applying (\ref{thm41b}) in the proof of Theorem \ref{thm41}, we get
by (\ref{exa52a})
$$
B\le \sqrt{ 2\int_\nn (\lambda +n)\tilde h^2 \mu(dn)} \sqrt{I} \le
\sqrt{2 [(\sqrt{\lambda}+1)^2+\lambda]} \ I.
$$
Thus we have proven

\bcor For the $M/M/\infty$ queue, if the Lipschitzian norm
$\|g\|_{Lip}$ of $g$ w.r.t. the Euclidean metric is finite (and
$\mu(g)=0$), then (\ref{key1}) and Bernstein's inequality
(\ref{Bern1}) hold with
$$
M=\|g\|_{Lip}\sqrt{2 \left[(\sqrt{\lambda}+1)^2+\lambda\right]}.
$$
\ncor
 }\nexa


\section{The subgeometric case}

\subsection{General result} In this last section, we will suppose no more that a
Poincar\'e inequality holds, and inspired by the Lyapunov function
approach, we introduce a more classical version of Lyapunov
condition
\begin{itemize}
\item[$(H_{LC})$] there exist a continuous function $U:\mathcal{X}\to
[1,+\infty)$ in $\dd_{e,\mu}(\LL)$, a measurable positive function
$\phi$,
 a set $C\in\BB$ with $\mu(C)>0$ and  constant $b>0$ such that
$$- \frac{\mathcal{L}U}{U}\ge \phi - b1_C,\  \mu\textrm{-a.s.}$$
\end{itemize}
In our mind $\phi$ goes to $0$ at infinity in this section.

We will also assume that a local Poincar\'e inequality holds for the
 set $C$ in ($H_{LC}$): there exists some constant $\kappa_C$ such
that for all $g\in \dd(\EE)$ such that $\mu(g1_C)=0$ \beq
\label{locSG} \mu (g^2 1_C)\le\kappa_C\mathcal{E}(g,g). \neq Note
that for diffusions on $\rr^d$, $C$ is often a ball $B(0,R)$ and the
local Poincar\'e inequality may then be easily deduced from the
local Poincar\'e inequality for the Lebesgue measure on balls by a
perturbation argument.

\bthm \label{thm6.1} Assume the Lyapunov function condition
$(H_{LC})$ and the local Poincar\'e inequality (\ref{locSG}) for the
set $C$. For $g\in L^2_0(\mu)$ such that $\sigma^2(g)$ is finite, if
$K_\phi(g^+)<+\infty$, then the transportation-information
inequality (\ref{key1}) holds with \beq\label{thm32alyap} M =
K_\phi(g^+)\left(b\kappa_C +1\right).\neq In particular Bernstein's
inequality (\ref{Bern1}) holds with that constant $M$. \nthm

\bprf In fact we have to slightly modify the key approach described
in section 2: for a constant $c>0$ to be chosen later,

\beq \label{key4}\aligned \nu(g)&=\int_\XX g h^2
d\mu = \int_\XX g
([h-c]^2 + 2 c h) d\mu\\
&=2c\<g, h\>_\mu + \int_\XX g[h-c]^2 d\mu=: A + B.
\endaligned
\neq For the first term $A=2c\<g, h\>_\mu$, since
$\sigma^2=\sigma^2(g)$ is assumed to be finite, we have by Remark
\ref{rem22} that $|A|\le c\sqrt{2\sigma^2 I}$.

Let consider now the second term
$$
B=\int_\XX g[h-c]^2 d\mu\le \int_\XX g^+[h-c]^2 d\mu\le K_\phi(g^+) \int_\XX
\left(b1_C- \frac {\LL U}{U} \right) [h-c]^2 d\mu.
$$
By a result in large deviations \cite[Lemma 5.6]{GLWY}, we have
$$
\int_\XX -\frac {\LL U}{U}[h-c]^2 d\mu\le \EE(h,h)=I.
$$
For the other term we apply the local Poincar\'e inequality, valid if we consider $c=\mu(h1_C)$ which leads to
$$
B\le K_\phi(g^+)\left(b\kappa_C +1\right) I.
$$
Remark finally that $c=\mu(h1_C)\le 1$. \nprf

Now we present an easy sufficient condition for the finiteness of
$\sigma^2(g)$ (and then for the CLT by Remark \ref{rem22}) by
following Glynn and Meyn \cite{GMeyn96}, which has its own interest.

\blem\label{lem62} Suppose that $R_1=\int_0^\infty e^{-t}P_t dt$ is
$\mu$-irreducible (i.e. $\mu\ll R_1(x,\cdot)$ for every $x\in \XX$)
and Harris positive recurrent (\cite{MT}). Assume that there are

$\bullet$ a (Lyapunov) continuous function $W : \XX\to [1,+\infty)$
in the extended domain $\dd_{e}(\LL)$ (see \S 4.1),

$\bullet$ a measurable function $F:\XX\to (0,+\infty)$,

$\bullet$ a $R_1$-small set $C$ with $\mu(C)>0$, i.e. $R_1(x,A)\ge
\delta \nu(A)$ for all $x\in C, A\in\BB$ for some constant
$\delta>0$ and $\nu\in \MM_1(\XX)$,

$\bullet $ and a positive constant $b$

\noindent such that $W$ is bounded on $C$ and \beq\label{lem62a} \LL
W\le -F + b 1_{C}. \neq If $|g|\le cF$ for some constant $c>0$ and
$\mu(g)=0$, then

\begin{enumerate}
\item There exists some measurable function $G$ such that $|G|\le
cW$ for some constant $c>0$, such that for any $t>0$, $\int_0^t
P_s|g| ds <+\infty$ and $P_tG-G=-\int_0^t P_sg ds$ everywhere on
$\XX$ (in such case we say that $G$ belongs to the extended domain
in the strong sense $\dd_s(\LL)$ of $\LL$ and write $-\LL G=g$).

\item If furthermore $g\in L^p_0(\mu)$ and $W\in L^q(\mu)$ where $p\in [2,+\infty]$ and $1/p+1/q=1$,
then $\sigma^2(g)$ is finite. \end{enumerate} \nlem

Its proof is postponed to the Appendix.

\subsection{Particular case: diffusions on $\rr^d$}
We study here the diffusion in $\rr^d$ with generator $\LL=\Delta
-\nabla V\cdot\nabla$ and $\mu=e^{-V} dx/Z$, presented in Section 4.
The first thing to remark is that any compact set is a small set,
and thus balls are small sets. A local Poincar\'e inequality is then
available. We then have

\bcor\label{cor63} Suppose that there exists a positive and bounded
function $\tilde \phi$ such that
\begin{equation}\label{kustr2}
\exists a<1, R, c>0, {\rm such~that~if~ }\ |x|>R,\qquad (1-a)|\nabla
V|^2-\Delta V\ge \tilde \phi(x).
\end{equation}
Then the weak Lyapunov condition $(H_{LC})$ is satisfied with
$U=e^{aV}$ with $\phi= a \tilde \phi$ and $C=B(0,R)$; and if $\int
e^{(a-1)V} dx<+\infty$ (i.e. $\mu(U)<+\infty$), then for any $\mu$
centered bounded function $g$ such that $|g|\le c_1 \tilde \phi U$
and $g(x)\le c_2 \tilde \phi$ for some positive constants $c_1,c_2$,
the asymptotic variance $\sigma^2(g)$ is finite by Lemma \ref{lem62}
and Bernstein's inequality holds. \ncor Note that, in parallel to
the second condition of Corollary \ref{cor1}, one may also consider
Lyapunov function of the form $U(|x|)$, but the result is then not
as explicit and we prefer to illustrate such an approach through
examples.

\bexa{\rm {\bf (sub-exponential measure)} Let $V(x)=|x|^\beta$ (if
$|x|>1$) for $\beta\in (0,1)$  such that no Poincar\'e inequality
holds. However, one may apply the previous corollary with
$U(x)=e^{a|x|^\beta}$ and  $\tilde\phi(x)=(1-a-\delta) \beta^2
(1+|x|)^{2(\beta-1)}$
 ($a,\delta \in(0,1)$, $a+\delta<1$). Hence by Corollary
 \ref{cor63},
Bernstein's inequality holds for $\mu$ centered bounded function $g$
such that for large $|x|$, $g(x)\le c/(1+|x|)^{2(1-\beta)}$. }\nexa

\bexa{\rm {\bf (Cauchy type measure)} Let $V(x)=\frac 12
(d+\beta)\log(1+|x|^2)$ for $\beta>0$. The condition $(H_{LC})$
holds with $U=e^{aV}= (1+|x|^2)^{a(d+\beta)/2}$ and
$\tilde\phi(x)=c/(1+|x|^2)$ for some constant $c>0$, where $a\in
(0,1)$ so that $(1-a)(d+\beta)>d$ (for $\mu(U)<+\infty$). So
Bernstein's inequality holds for $\mu$ centered bounded function $g$
such that for large $|x|$, $g(x)\le K/(1+|x|^{2})$ for some constant
$K>0$, by Corollary
 \ref{cor63}. }\nexa

\brmk{\rm  One may be surprised that the upper bound for the test
function is the same for every Cauchy type measure. One may find the
beginning of an answer in recent results of Bobkov-Ledoux \cite{BL}
(see also Cattiaux-Gozlan-Guillin-Roberto \cite{CGGR}).  Indeed, in
their work they prove that this type of measures satisfy a weighted
Poincar\'e type inequality where the weight is the same for every
Cauchy-type measure.} \nrmk

\subsection{Particular case : birth-death processes}
We adopt here the notations of subsection 4.3, and assume once again
that the process is positive recurrent. We suppose for simplicity
that for large enough $n$, the death rate $a_n$ is larger than the
birth rate $b_n$.

\bcor\label{cor65} If there are $m> 0$, $N\ge 1$ and a positive
sequence $(c_n)_{n\in\nn}$ such that
\begin{enumerate}
\item for all $n\ge N$, $a_n-b_n\ge c_n>0$;

\item $\sum_n n^m \mu_n<+\infty$,
\end{enumerate}
then Bernstein's inequality is valid for every $\mu$ centered
bounded function $g$ such that  for large $n$, $|g(n)|\le c
n^{m-1}c_n$ and $g(n)\le Kc_n/n$ for some constants $c,K>0$. \ncor
\bprf{\rm Let $U(n)=(1 + n)^m$, then for large $n$,
$$-\frac{\LL U(n)}{U(n)}\ge m (a_n-b_n) \left(\frac 1n + o(\frac 1n)\right).$$
Hence the Lyapunov condition ($H_{LC}$) holds for
$\phi(n)=(m-\delta) c_n/(1+n)$ where $\delta\in (0,m)$. The local
Poincar\'e inequality is always valid in this context and a precise
estimation of the constant may be found in Chen \cite{Chen05}. Since
$\mu(U)$ is finite, we can apply Lemma \ref{lem62} to conclude that
$\sigma^2(g)<+\infty$ for $|g|\le c \phi U$. It remains to apply
Theorem \ref{thm6.1}.

}\nprf

\bexa {\rm Let $b_n\equiv 1$ and $a_n = 1+ a/(n+1)$ where $a>0$.
Then $c_n:=a_n-b_n=a/(n+1)$ and $\pi_n$ behaves as $\frac 1{n^a}$
for large $n$. Thus the process is positive recurrent if and only if
$a>1$. For $a>1$, take $m\in (0, a-1)$, we see that the conditions
in Corollary \ref{cor65} are all satisfied. Hence Bernstein's
inequality holds for $\mu$-centered $g$ such that $|g(n)|\le K/n^2$
for large $n$. This is quite similar as in the Cauchy measure case.

 } \nexa

\section{Appendix}

\bprf[Proof of Lemma \ref{lem62}] Let us first prove part (2) by
admitting part (1). Let $G$ be the strong solution of $-\LL G=g$
given in part (1). Since $W\in L^q(\mu)$, considering $G-\mu(G)$ if
necessary we may and will assume that $\mu(G)=0$. Now for any
$\vep>0$, let $R_\vep=\int_0^\infty e^{-\vep t} P_t
dt=(\vep-\LL)^{-1}$ be the resolvent. By the resolvent equation, $G-
R_\vep g=\vep R_\vep G$ which tends to $\mu(G)=0$ in $L^q(\mu)$ as
$\vep\to 0$ by the ergodic theorem, we have
$$
\lim_{\vep\to 0}\<R_\vep g,g\>_\mu = \int G g d\mu <+\infty.
$$
This relation yields that $\sigma^2(g)$ in (\ref{variance}) exists
and $\sigma^2(g)= 2  \int G g d\mu$ (in the actual symmetric case).

We turn now to prove part (1). This is due to
 Glynn and Meyn
\cite[Theorem 3.2]{GMeyn96} when $F$ is bounded from below by a
positive constant. Let us modify slightly their proof for the
general case.

{\bf Step 1 (Reduction to the discrete time case).} At first since
$e^{-bt}W(X_t)$ is a local super-martingale, then a
super-martingale, so $P_tW\le e^{bt}W$ for all $t\ge0$. Moreover for
any $\lambda>0$, by It\^o's formula,
$$
M_t=e^{-\lambda t} W(X_t) - W(X_0)+\int_0^t e^{-\lambda s}
\left(\lambda W -\LL W\right)(X_s)ds
$$
is a $\pp_x$-local martingale for every $x\in \XX$. Hence taking a
sequence of stopping times $(\tau_n)$ increasing to $+\infty$ such
that $\ee^xM_{\tau_n}=0$, we have for every $x\in\XX$,
$$
\ee^{x} \int_0^{\tau_n}  e^{-\lambda s} \left( \lambda W + F - b1_C
\right)(X_s) ds \le \ee^x \int_0^{\tau_n} e^{-\lambda s}
\left(\lambda W -\LL W\right)(X_s)ds\le W(x).
$$
Letting $n$ go to infinity, we obtain by monotone convergence
$$\lambda R_\lambda W + R_\lambda F\le W
 + b R_\lambda 1_{C}.$$

Consider the Markov kernel $Q=R_1$. The relation above says that
\beq\label{lem62b} Q W \le W -Q F + b Q 1_C. \neq

Assume that one can prove that there is $G$ such that $|G|\le c W$
(for some constant $c>0$) such that \beq\label{lem62c} (1-Q) G =
Qg.\neq Then $G=R_1(G+g)\in \dd_s(\LL)$ and
$R_1(-\LL)G=(1-R_1)G=R_1g$. Consequently $-\LL G= (1-\LL)R_1(-\LL)G=
(I-\LL)R_1g=g$, the desired claim in part (1).

Therefore it remains to solve (\ref{lem62c}) under the condition
(\ref{lem62b}).

{\bf Step 2 (atom case).} Let us suppose at first that the small set
$C$ in (\ref{lem62b})  is an atom of $Q$, i.e.,
$Q(x,\cdot)=Q(y,\cdot)$ for all $x,y\in C$. In this case one
solution to (\ref{lem62c}) is given by \beq\label{lem62d} G(x) =
\ee^x \sum_{k=0}^{\sigma_C} Qg(Y_k) \neq where $(Y_n)_{n\ge0}$ is
the Markov chain with transition probability kernel $Q$ defined on
$(\Omega, (\FF_n), \qq_x)$ equipped with the shift $\theta$ (so that
$Y_n(\theta \omega)=Y_{n+1}(\omega)$), $\sigma_C=\inf\{n\ge0; \
Y_n\in C\}$.

To justify this fact which is one key in \cite{GMeyn96}, notice

1) $G$ given  by (\ref{lem62d}) is well defined. In fact $|Qg|\le c
QF$. Using the condition (\ref{lem62b}) and the fact that
$$
W(Y_n) - W(Y_0) + \sum_{k=0}^{n-1} (W-QW)(Y_k)
$$
is a $\qq_x$-martingale, we obtain the following at first for
$\sigma_C\wedge n$ and then for $\sigma_C$ (by letting $n\to\infty$)
 $$\aligned\ee^x
\sum_{0\le k\le \sigma_C -1 } QF(Y_k)&\le b\ee^x \sum_{0\le k\le
\sigma_C -1 } Q1_C(Y_k) + W(x)\\
&= b\ee^x \sum_{1\le k\le \sigma_C } 1_C(Y_k) +W(x)\le b+W
\endaligned$$
where the second equality for $\sigma_C\wedge n$ (instead of
$\sigma_C$) follows by Doob's stopping time theorem. Consequently
$$\ee^x \sum_{0\le k\le \sigma_C } QF(Y_k)\le \sup_{x\in C} QF(x) +
\ee^x \sum_{k=0}^{\sigma_C-1} QF(Y_k)\le \sup_{x\in C} QF(x) + b +
W(x).$$ By (\ref{lem62b}), $QF\le W +b$ is bounded on $C$. Therefore
$G$ is well defined and $|G|\le c(b' + W)$.

2) Let $\tau_C:=\inf\{n\ge 1; Y_n \in C\}$. We have
$\sigma_C\circ\theta =\tau_C-1$ on $[\sigma_C=0]$ and
$\sigma_C\circ\theta =\sigma_C-1$ on $[\sigma_C\ge 1]$. Hence for
$x\in C$
$$
QG(x) = \ee^x \sum_{k=0}^{\sigma_C\circ \theta} Qg(Y_{k+1})= \ee^x
\sum_{k=1}^{\tau_C} Qg(Y_{k})
$$
which is constant on $x\in C$ and equals to $\mu(g)/\mu(C)=0$, then
$G(x)-QG(x)=G(x)=Qg(x)$ for $x\in C$. Now for $x\notin C$,
$$
QG(x) = \ee^x \sum_{k=0}^{\sigma_C\circ \theta} Qg(Y_{k+1})= \ee^x
\sum_{k=0}^{\sigma_C-1} Qg(Y_{k+1})= G(x)-Qg(x).
$$
So $G-QG=Qg$ everywhere on $\XX$.

{\bf Step 3 (non-atom case).} In the non-atom case one can consider
the splitting chain in \cite[Proof of Theorem 2.3]{GMeyn96} to
reduce the problem to the atom case.
 \nprf

\end{document}